\input amstex
\magnification =\magstep 1
\documentstyle{amsppt}
\pageheight{9truein}
\pagewidth{6.5truein}
\NoRunningHeads
\baselineskip=16pt

\topmatter
\title Another criterion for solvability of    finite groups
\endtitle

\author Marcel Herzog*, Patrizia Longobardi** and Mercede Maj**
\endauthor

\affil *School of Mathematical Sciences \\
       Tel-Aviv University \\
       Ramat-Aviv, Tel-Aviv, Israel
{}\\
       **Dipartimento di Matematica \\
       Universit\`a di Salerno\\
       via Giovanni Paolo II, 132, 84084 Fisciano (Salerno), Italy
\endaffil

\thanks This work was supported by the National Group for Algebraic and
Geometric Structures, and their Applications (GNSAGA - INDAM), Italy.
\endthanks

\date 19/09/2021
\enddate

\abstract Let $G$ be a finite group. Denote by $\psi(G)$ the sum
$$\psi(G)=\sum_{x\in G}|x|,$$
where $|x|$ denotes the order of the element $x$, and
by $o(G)$ the quotient
$$o(G)=\frac{\psi(G)}{|G|}.$$
Confirming a conjecture posed by
E.I. Khukhro, A. Moreto and M. Zarrin, we prove that
if $o(G)< o(A_5)$, then  $G$ is solvable.
\endabstract

\endtopmatter

\document

\heading 1. Introduction.\\
\endheading

Let $G$ be a finite group. Denote by $\psi(G)$ the sum
$$\psi(G)=\sum_{x\in G}|x|,$$
where $|x|$ denotes the order of the element $x$, and by $o(G)$ the quotient
$$o(G)=\frac{\psi(G)}{|G|}.$$
Thus $o(G)$ denotes the average element order of $G$.
Moreover, if $S\subseteq G$, then we define $\psi(S)=\sum_{x\in S}|x|.$

Recently many authors studied the function $\psi(G)$ and, more generally, properties of finite groups
determined by their element orders (see for example [1]-[5], [9], [12]-[17], [20], [22]-[29], [34], [36] and [38]).

In their paper [20], E.I. Khukhro, A. Moreto and M. Zarrin posed the following conjecture:
\proclaim{Conjecture} Let $G$ be a finite group and suppose that
$$o(G)<o(A_5).$$
Then $G$ is solvable.
\endproclaim

Notice that
$$o(A_5)=\frac {\psi(A_5)}{|A_5|}=\frac {211}{60}=3.51666... $$

In this paper we prove that their conjecture is true. In fact, we prove the following three theorems.
\proclaim{Theorem A}Let $G$ be a finite group and suppose that
$$o(G)<o(S_3)=\frac {13}{6}.$$
Then $G$ is an elementary abelian $2$-group.

Moreover, if $o(G)=\frac {13}6$, then $G\simeq S_3$.
\endproclaim.
\proclaim{Theorem B} Let $G$ be a finite group and suppose that
$$o(G)\leq o(A_5)=\frac {211}{60}.$$
Then either $G\simeq A_5$, or $G$ is solvable.
\endproclaim

Theorem B immediately implies the main result of this paper.
\proclaim{Theorem C} Let $G$ be a finite group and suppose that
$$o(G)< o(A_5)=\frac {211}{60}.$$
Then  $G$ is solvable.
\endproclaim

We shall use the following notation. If $G$ is a finite group, then $1$ will denote the identity  element of $G$ and
sometimes also the group $\{1\}$.
We shall denote by $i_2(G)$ the number of
elements of $G$ of order $2$ and by $i_3(G)$ the number of
elements of $G$ of order  $3$. Sometimes we shall use the shorter notation $i_2$ and $i_3$, if there is no ambiguity.
We shall denote the set of primes dividing $|G|$ by $\Pi(G)$ and the set of primes dividing a positive integer $n$ by $\pi(n)$.
Moreover, simple groups will be assumed to be non-abelian.
The rest of our  notation is the usual one (see for example [10], [19] and [37]).

In Section 2 we shall  recall some known results concerning $i_2(G)$ and $i_3(G)$.

In Section 3 we shall prove the following key lemma.
\proclaim{Lemma 3.3} Let $G$ be a non-solvable finite group and let $p$ be a prime dividing $|G|$.
Then the following statements hold.
\roster
\item If $p\geq 17$ and $o(G)\leq o(A_5)$, then $G$ is $p$-solvable.
\item If $p\geq 23$ and $o(G)<3.55$, then $G$ is $p$-solvable.
\item If $p\geq 11$ and $o(G)<3.4479$, then $G$ is $p$-solvable.
\endroster
\endproclaim

The proof of Theorem B will be achieved in two steps. The first step is to study {\it simple} groups $G$ satisfying $o(G)\leq o(A_5)$.
By Lemma 3.3(1), such simple groups  satisfy the following important condition:
$$\Pi(G)\subseteq \{2,3,5,7,11, 13\}.$$
Simple groups satisfying this condition will be called {\it our simple groups}. These groups belong
to one of the following four classes: $ |\Pi(G)|=3,4,5\ \text{or}\ 6$.
Simple  groups which belong
to each such class have been studied by various authors, as indicated in the preface to Proposition 2.5. Using these studies
and taking into account  that $\Pi(G)\subseteq \{2,3,5,7,11, 13\}$, we present in Proposition 2.5 a list of
simple groups  which contains all our simple groups. This list will be referred to as {\it The List}.

Our aim in Section 4 is to prove the following theorem.
\proclaim {Theorem 4.19} Let $G$ be a finite simple group with $\Pi(G)\subseteq \{2,3,5,7,11,13\}$.
Then either $G\simeq A_5$ or $o(G)\geq 3.55> o(A_5)$.
\endproclaim

So the simple groups $G\not\simeq A_5$  which belong to The List do not satisfy $o(G)\leq o(A_5)$. Therefore, in order to prove Theorem B, it suffices to show
that there does not exist  a non-simple  non-solvable group $G$ satisfying  $o(G)\leq o(A_5)$. This is accomplished in Section 6,
where the proof of Theorem B is completed. Also Theorem A is proved in Section 6.

For the proofs of Theorems A and B we need a bound on the proportion of elements in a  non-solvable  group inverted by an automorphism.
In Section 5, generalizing a previous result of W.M. Potter (see [35], Theorem 3.3), while using his methods and his ideas, we prove the following theorem.
\proclaim {Theorem 5.5} Let $G$ be a finite group which contains no non-trivial normal solvable subgroups and let
$\theta \in Aut(G)$. If $\theta$ inverts more than $\frac 29|G|$ elements of $G$, then $G\simeq A_5$.
\endproclaim

Finally we prove some basic results concerning the function $o(G)$.
\proclaim {Lemma 1.1} Let $G$ be a finite group and $G\neq 1$. Then the following statements hold.

\noindent (1) We have $o(G)\geq 2-\frac 1{|G|}\geq \frac 32$. In particular, if $G$ is an elementary abelian $2$-group, then
$o(G)=2-\frac 1{|G|}$ and if $G$ is not an elementary abelian $2$-group, then
$o(G)\geq 2+\frac 1{|G|}$. Hence $o(G)\leq 2$ if and only if $G$ is an elementary abelian $2$-group and $o(G)=2-\frac 1{|G|}$.

\noindent (2) If $G$ is of odd order, then
$o(G)\geq 3-\frac 2{|G|}\geq 3-\frac 23=\frac 73$.

\noindent (3) If $G=A\times B$ with $(|A|,|B|)=1$, then $o(G)=o(A)o(B)$. In particular, if $A\neq 1$  and $B\neq 1$,
then $o(G)\geq \frac 72$.
\endproclaim

\demo{Proof}(1) Clearly $\psi(G)\geq 2(|G|-1)+1=2|G|-1$. Hence $o(G)\geq 2-\frac 1{|G|}\geq  \frac 32$,
and if $G$ is an elementary abelian $2$-group, then $o(G)=2-\frac 1{|G|}$.
If $G$ is not an elementary abelian $2$-group, then there exists $x\in G$ with $|x|=|x^{-1}|\geq 3$. Since $x\neq x^{-1}$,
it follows that $\psi(G)\geq (2|G|-1)+2\cdot1=2|G|+1$ and hence $o(G)\geq 2+\frac 1{|G|}>2$. Therefore
$o(G)\leq 2$ if an only if $G$ is an elementary abelian $2$-group and $o(G)=2-\frac 1{|G|}$.

\noindent
(2) If $G$ is of odd order, then $\psi(G)\geq 3(|G|-1)+1=3|G|-2$ and $o(G)\geq 3 -\frac 2{|G|}\geq 3-\frac 23=\frac 73$.

\noindent (3)  If $G=A\times B$ with $(|A|,|B|)=1$, then  $|G|=|A||B|$ and
$\psi(G)=\psi(A)\psi(B)$. Hence $o(G)=o(A)o(B)$.
If $|A|\neq 1$ and $|B|\neq 1$, then either $A$ or $B$ is of odd order, and it follows by (1) and (2) that $o(G)\geq \frac 32\cdot \frac 73=\frac 72.$
\qed
\enddemo

We shall use the results of Lemma 1.1 freely throughout this paper.
\heading 2. Some known results.\\
\endheading
Let $G$ be a finite group. Recall that  we denote by $i_2(G)$ the number of elements of $G$ of order $2$ and by $i_3(G)$
the number of elements of $G$ of order $3$.

Let
$$T(G)=\sum_{\chi \in Irr(G)}\chi(1).$$
It is well known that
$$i_2(G)+1\leq T(G)$$
(see for example [21] and [33]).

We shall need the following related result (see  Theorem A in [18]).

\proclaim{Proposition 2.1} Let $p\geq 7$ be a prime. Write
$$g_1(p)=p-1\quad \text{if}\quad p\equiv 3\pmod 4,$$
$$g_1(p)=\frac {p(p^2-1)}{p^2+p+2}\quad \text{if}\quad  p\equiv 1\pmod 4,$$
and
$$g_2(p)=\frac {p(p-2)}{p-1}.$$
If $G$ is a finite group, $p\geq 7$ is a prime and $G$ is non-$p$-solvable, then one of the
following cases holds:
$$T(G)=\frac {|G|}{g_1(p)},\quad T(G)=\frac {|G|}{g_2(p)},\quad  T(G)\leq \frac {|G|}{p-1}.$$
\endproclaim

Proposition 2.1 yields the following corollary.

\proclaim{Corollary 2.2} Let $p$ be a prime and let $G$ be a  non-$p$-solvable finite group.
Then the following statements hold.
\roster
\item If $p\geq 11$, then $i_2(G)<\frac {10}{99}|G|$.
\item If $p\geq 17$, then $i_2(G)<\frac 1{15}|G|$.
\item If $p\geq 23$, then $i_2(G)<\frac 1{20}|G|$.
\endroster
\endproclaim
\demo{Proof}  Notice that $g_1(p),g_2(p)>p-2$ and
by the previous remarks
$$i_2(G)< T(G).\tag{$*$}$$
Since in all cases $p\geq 11$,  Proposition 2.1 implies that
$$i_2(G)<T(G)<\frac 1{p-2}|G|.\tag{$**$}$$

 Let $p=11$. Since $g_1(11)=10$ and $g_2(11)=\frac {11(11-2)}{10}=\frac {99}{10}$, it follows by
Proposition 2.1 that  $T(G)\leq \frac {10}{99}|G|$ and hence by ($*$) $i_2(G)<\frac {10}{99}|G|$.
If $p\geq 13$, then by ($**$) $i_2(G)<\frac 1{11}|G|<\frac {10}{99}|G|$. This proves (1).

The  inequalities in (2)-(3) hold by ($**$).
\qed
\enddemo

We also need the following two results.
\proclaim{Proposition 2.3} Let $G$ be a  finite group and let
$\phi$ be an automorphism of $G$. Then the following statements hold.
\roster

\item If $\phi$ inverts more than $\frac 4{15}|G|$ elements of $G$, then $G$ is solvable.
\item If $\phi$ inverts more than $\frac 34|G|$ elements of $G$, then $G$ is abelian.
\endroster
\endproclaim

\demo{Proof} See [35], Corollary 3.2 and Corollary 2.4.
\enddemo

\proclaim{Proposition 2.4} If $G$ is a non-solvable  finite group. then
$$i_2(G)\leq \frac 4{15}|G|-1\quad \text{and}\quad i_3\leq \frac 7{20}|G|-1.$$
\endproclaim
\demo{Proof} See [7], Lemma 2.9 and Lemma 2.16.
\enddemo

In Section 4, we shall deal with our simple groups. Recall that a simple group $G$ is called an our simple group
if
 $\Pi(G)\subseteq \{2,3,5,7,11,13\}$, so $3\leq |\Pi(G)|\leq 6$.
Finite simple groups with $|\Pi(G)|=3$ have been studied by M. Herzog in [11], finite simple groups with $|\Pi(G)|=4$
have been studied by Y. Bugeaud, Z. Cao and M. Mignotte in [6], and finite simple groups with $|\Pi(G)|=5,6$
have been studied by A. Jafarzadeh and A. Iranmanesh in [30]. From their results, using also the information in the
paper [31] of D. Yu, J. Li, G. Chen, L. Zhang and W. Shi, namely Lemma 2.1 for the case $|\Pi(G)|=3$, Lemmas 2.2-2.4 for the case  $|\Pi(G)|=4$,
Lemmas 2.6-2.10 for the case $|\Pi(G)|=5$, and in the Atlas [8], pages 239-241 (aided by some computations)
for the case $|\Pi(G)|=6$, the following proposition
follows. Recall that if $n$ is a positive integer, then $\pi(n)$ denotes the set of primes dividing $n$.

\proclaim{Proposition 2.5} Let $G$ be a finite simple group with $\Pi(G)\subseteq \{2,3,5,7,11,13\}$.
Then the following statements hold.
\roster
\item If $|\Pi(G)|=3$, then $G$ is isomorphic to one of the following groups:
$$A_5,A_6,S_4(3)\simeq U_4(2),L_2(7),L_2(8),L_3(3), U_3(3).$$
\item If $|\Pi(G)|=4$, then $G$ is isomorphic to one of the following groups:
$$ \align
L_2(q),\ &\text{q a prime},\ |\pi(q^2-1)|=3,\\
L_2(2^m),\ &2^m-1=u,\ 2^m+1=3t,\ \text{u,t are primes},\ t>3,\\
L_2(3^m),\ &3^m-1=2u,\ 3^m+1=4t,\ \text{u,t are primes},\\
L_2(25),\ &L_2(49),\ L_3(4),\ L_4(3),\ U_3(4),\ U_3(5),\ U_4(3),\ U_5(2),\\
S_4(5),\ &
S_4(7),\ S_6(2),\ O_8^+(2),\ G_2(3),\ ^3D_4(2),\ ^2F_4(2)',\ Sz(8),\\
A_7,\ &A_8,\ A_9,\ A_{10},\ M_{11},\ M_{12},\ J_2.
\endalign
$$
\item If $|\Pi(G)|=5$, then $G$ is isomorphic to one of the following groups:
$$ \align
L_2(q),\ &\text{q a prime power},\ |\pi(q^2-1)|=4,\\
L_3(q),\ &\text{q a prime power },\ |\pi((q^2-1)(q^3-1))|=4,\\
U_3(q),\ &\text{q a prime power },\ |\pi((q^2-1)(q^3+1))|=4,\\
O_5(q)&\simeq S_4(q),\ \text{q a prime power },\ |\pi(q^4-1)|=4,\\
Sz(2^{2m+1})&\simeq \ ^2B_2(2^{2m+1}),\ \text{with}\ |\pi((2^{2m+1}-1)(2^{4m+2}+1))|=4,\\
R(q),\ &\text{with}\ q=3^{2m+1},\ |\pi(q^2-1)|=3\ \text{and}\ |\pi(q^2-q+1)|=1,\\
L_5(3),\ &S_6(3),\ U_4(5),\ U_6(2),\ O_7(3),\ O_8^+(3),\ G_2(4),\\
A_{11},\ &A_{12},\ M_{22},\ HS,\ McL.
\endalign
$$
\item If $|\Pi(G)|=6$, then $G$ is isomorphic to one of the following groups:
$$ \align
L_2(q),\ &\text{q a prime power},\ |\pi(q^2-1)|=5,\\
L_3(q),\ &\text{q a prime power },\ |\pi((q^2-1)(q^3-1))|=5,\\
L_4(q),\ &\text{q a prime power },\ |\pi((q^2-1)(q^3-1)(q^4-1))|=5,\\
U_3(q),\ &\text{q a prime power },\ |\pi((q^2-1)(q^3+1))|=5,\\
U_4(q),\ &\text{q a prime power },\ |\pi((q^2-1)(q^3+1)(q^4-1))|=5,\\
O_5(q)&\simeq S_4(q),\ \text{q a prime power },\ |\pi(q^4-1)|=5,\\
G_2(q),\ &\text{q a prime power },\ |\pi(q^6-1)|=5,\\
Sz(2^{2m+1})&\simeq \ ^2B_2(2^{2m+1}),\ \text{with}\ |\pi((2^{2m+1}-1)(2^{4m+2}+1))|=5,\\
R(3^{2m+1})&,\ \text{with}\ |\pi((3^{2m+1}-1)(3^{6m+3}+1))|=5,\\
L_6(3),\ &A_{13},\ A_{14},\ A_{15},\ A_{16,},\ Suz,\ Fi_{22}.
\endalign
$$
\endroster
\endproclaim

\heading 3. Some preliminary results.\\
\endheading

In this section we shall prove three useful lemmas, which we shall use
in our  proofs of Theorems A and B.

\proclaim{Lemma 3.1} Let $G$ be a finite group containing a non-trivial normal subgroup $H$.
Then the following statements hold.
\roster
\item If $x\in G\setminus H$, then the order  $|xH|$ of $xH$ in $G/H$
divides the order of $xh$ in $G$ for every $h\in H$. In particular, $|xh|\geq |xH|$
for every $h\in H$.
\item $o(G/H)<o(G)$ .
\endroster
\endproclaim
\demo{Proof} (1) If $h\in H$ and $|xh|=n$, then $(xh)^n=1$ and $(xH)^n=(xhH)^n=H$.
Hence $|xH|$ in $G/H$ divides $n$, as claimed. In particular, $|xh|\geq |xH|$
for every $h\in H$.

(2) Write $G/H=\{H,x_1H,\dots,x_sH\}$. Then
$$G=H\dot\cup x_1H \dot\cup\dots\dot\cup x_sH\qquad \text{and}\qquad
\psi(G)=\psi(H)+\psi(x_1H)+\dotsb +\psi(x_sH).$$
By (1), $\psi(x_jH)\geq |H||x_jH|$ for every $j\in \{1,\dots,s\}$. Therefore
$$\psi(G)\geq \psi(H)-|H|+|H|(1+|x_1H|+\dotsb +|x_sH|)=\psi(H)-|H|+|H|\psi(G/H),$$
and
$$o(G)\geq \frac {\psi(H)-|H|}{|G|}+\frac {|H|}{|G|}\psi(G/H)=
\frac {\psi(H)-|H|}{|G|}+o(G/H)>o(G/H),$$
since $\psi(H)>|H|$.
\qed
\enddemo

\proclaim{Lemma 3.2} Let $G$ be a non-solvable finite group. Then the following statements hold.
\roster
\item $o(G)>3.11$.
\item If $i_2(G)\leq \frac {|G|}{20}$, then $o(G)\geq 3.55$.
\item If $i_2(G)\leq \frac {10}{99}|G|$, then $o(G)>3.4479$.
\item If $i_2(G)< \frac {|G|}{15}$, then $o(G)> o(A_5)$.
\item If $i_2(G)\leq \frac {|G|}{16}$ and $i_3(G)+1\leq \frac {|G|}{14}$, then $o(G)>3.8$.
\endroster
\endproclaim

\demo{Proof} Write $R=\{x\in G\mid |x|\geq 4\}$ and let $r=|R|$. Then
$$|G|=1+i_2(G)+i_{3}(G)+r,\quad \text{yielding}\quad r-i_2=|G|-2i_2-(i_3+1).$$
Since $G$ is non-solvable, there exist  elements $y$ and $y^{-1}$ in $R$ with $|y|\geq 5$.
Hence
$$\psi(G)\geq 1+2i_2+3i_3 +4r+2=3+3i_2+3i_3+3r +(r-i_2).$$
Thus
$$\psi(G)\geq 3|G|+(r-i_2).\tag{$*$}$$
Recall that since $G$ is non-solvable,  Proposition 2.4 implies that $i_2<\frac 4{15}|G|$ and $i_3+1\leq \frac 7{20}|G|$.

(1) We have
$$r-i_2=|G|-2i_2-(i_3+1)>|G|-(\frac 8{15}+\frac 7{20})|G|=|G|-\frac {53}{60}|G|=\frac 7{60}|G|> (0.11)|G|$$
and it follows by ($*$) that $\psi(G)> (3.11)|G|$. Hence $o(G)>3.11$, as required.

(2) If $i_2\leq \frac {|G|}{20}$, then
$$r-i_2=|G|-2i_2-(i_3+1)\geq |G|-(\frac 2{20}+\frac 7{20})|G|=|G|-\frac 9{20}|G|=\frac {11}{20}|G|=(0.55)|G|$$
and it follows by ($*$) that $\psi(G)\geq (3.55)|G|$. Hence $o(G)\geq 3.55$, as required.

(3) If $i_2\leq \frac {10}{99}|G|$, then
$$r-i_2=|G|-2i_2-(i_3+1)\geq |G|-(\frac {20}{99}+\frac 7{20})|G|=|G|-\frac {1093}{1980}|G|=\frac {887}{1980}|G|.$$
So $r-i_2>0.4479|G|$ and it follows by ($*$) that $\psi(G)> (3.4479)|G|$. Hence $o(G)> 3.4479$, as required.

(4) If $i_2< \frac {|G|}{15}$, then
$$r-i_2=|G|-2i_2-(i_3+1)>|G|-(\frac 2{15}+\frac 7{20})|G|=|G|-\frac {29}{60}|G|=\frac {31}{60}|G|$$
and it follows by ($*$) that $\psi(G)> (3+\frac {31}{60})|G|=\frac {211}{60}|G|=o(A_5)|G|$. Hence $o(G)>
 o(A_5)$, as required.

(5) In this case, we have
$$r-i_2=|G|-2i_2-(i_3+1)\geq |G|-(\frac 2{16}+\frac 1{14})|G|=|G|-\frac {22}{112}|G|=\frac {90}{112}|G|>(0.8)|G|$$
and it follows by ($*$) that $\psi(G)>(3.8)|G|$. Hence $o(G)> 3.8$, as required.
\qed
\enddemo

\proclaim{Lemma 3.3} Let $G$ be a non-solvable finite group and let $p$ be a prime dividing $|G|$.
Then the following statements hold.
\roster
\item If $p\geq 17$ and $o(G)\leq o(A_5)$, then $G$ is $p$-solvable.
\item If $p\geq 23$ and $o(G)<3.55$, then $G$ is $p$-solvable.
\item If $p\geq 11$ and $o(G)<3.4479$, then $G$ is $p$-solvable.
\endroster
\endproclaim

\demo{Proof} (1) If $p\geq 17$ and $G$ is non-$p$-solvable, then $i_2(G)<\frac 1{15}|G|$ by Corollary 2.2(2).
This implies, by Lemma 3.2(4), that $o(G)> o(A_5)$, a contradiction. Therefore $G$ is $p$-solvable, as required.

Similarly, in the other cases, if $G$ is non-$p$-solvable we reach a contradiction to our assumptions
by Lemma 3.2.
Therefore in all cases $G$ is $p$-solvable, as required.
\qed
\enddemo

\heading 4. Our simple groups.\\
\endheading

 In this section $G$ denotes a finite simple group. We shall prove that if $G$ is  "our simple group"
 (i.e. $\Pi(G)\subseteq \{2,3,5,7,11,13\}$), then either $G\simeq A_5$ or $o(G)\geq 3.55>o(A_5)$.  Simple groups $G$ which satisfy $\Pi(G)\subseteq \{2,3,5,7,11,13\}$
 are described in Proposition 2.5. We shall prove our claim separately for groups $G$ satisfying $|\Pi(G)|=3,4,5, \text {or}\ 6$, as described in
 Proposition 2.5. This will be accomplished by means of numerous lemmas, in which groups $G$ which are mentioned in Proposition 2.5,
 other than $A_5$, are shown to satisfy  $o(G)\geq 3.55$. This will lead to the conclusion that among our simple groups, only $G\simeq A_5$
 satisfies $o(G)\leq o(A_5)$.
 In our proofs we rely heavily on the paper [7] of T.C. Burness and S.D. Scott and on the Atlas [8].

\proclaim{Lemma 4.1} Let $G\simeq A_n$, with $n\geq 5$. Then either $G\simeq A_5$ or $o(G)\geq 3.55$.
\endproclaim
\demo{Proof} First, let $n=6$. Then $\psi(A_6)=1411$ and $|A_6|=360$, so $o(A_6)=\frac {1411}{360}>3.91$, as required.

Suppose, now that $n\geq 7$. By Lemma 3.3 in [7],  we have
$$\frac {i_2(A_n)}{|A_n|}\leq \frac 2{8(n-4)!}+\frac 2{2^4\cdot 23}=\frac 14\left(\frac 1{(n-4)!}+\frac 1{46}\right).$$
Our aim is to prove that $\frac {i_2(A_n)}{|A_n|}\leq \frac 1{20}$, since then it follows by Lemma 3.2(2) that
$o(G) \geq 3.55$, as required.
 Since $\frac 14(\frac 1{(n-4)!}+\frac 1{46})$ is clearly a decreasing
function of $n$, it suffices to show that  $\frac {i_2(A_7)}{|A_7|}\leq \frac 1{20}$ . And indeed,
$$\frac {i_2(A_7)}{|A_7|}\leq \frac 14\left(\frac 16+\frac 1{46}\right)=\frac {13}{276}<\frac 1{20},$$
as required.
\qed
\enddemo

\proclaim{Lemma 4.2} Let $G\simeq L_2(q)$, with $q\geq 4$ and $\Pi(G)\subseteq\{2,3,5,7, 11,13\}$. Then
$$\text{either}\quad G\simeq A_5\quad \text{or}\quad o(G)\geq 3.55.$$
\endproclaim

\demo{Proof} If $q=4,5$, then $G\simeq A_5$. If $q=7$, then $\psi(L_2(7))=715$ and $|L_2(7)|=168$, so
$o(L_2(7))=\frac {715}{168}\geq 4.2$. If $q=8$, then $\psi(L_2(8))=3319$ and $|L_2(8)|=504$, so
$o(L_2(8))=\frac {3319}{504}\geq 6.5$. If $q=9$, then $L_2(9)\simeq A_6$ and $o(L_2(9))\geq 3.55$ by Proposition 4.1.
If $q=11$, then $\psi(L_2(11))=3741$ and $|L_2(11)|=660$, so $o(L_2(11))=\frac {3741}{660}\geq 5.6$.
Finally, if $q=13$, then $\psi(L_2(13))=7281$ and $|L_2(13)|=1092$, so $o(L_2(13))=\frac {7281}{1092}\geq 6.6$.

Since $17$ divides $|L_2(16)|$ and $|L_2(17)|$, and $19$ divides  $|L_2(19)|$, it remains only to consider groups $G=L_2(q)$ with
$q\geq 23$. If $q$ is even, then $|G|=q(q^2-1)$ and  $i_2(G)=q^2-1=\frac {|G|}q$,
and if $q$ is odd, then $|G|=\frac 12q(q^2-1)$ and $i_2(G)\leq \frac 12q(q+1)=\frac {|G|}{q-1}$ (see Lemma 3.4 in [7]).
Since $q\geq 23$,  in both cases $i_2(G)\leq \frac {|G|}{q-1}\leq \frac {|G|}{22}$, which implies by Lemma 3.2(2)
that $o(G)\geq 3.55$, as required.
\qed
\enddemo

\proclaim{Lemma 4.3} Let $G\simeq L_3(q)$, with $q=3,4$ or $q\geq 8$ .
Then
$$o(G)\geq 3.55.$$
\endproclaim
\demo{Proof}
If $q\in \{4,8,9\}$, then by [8] $G\simeq L_3(q)$ has only one conjugacy class of involutions and the orders of their Sylow $2$-subgroups
are $\{2^6,2^9,2^7\}$, respectively.  If $q=3$, then by [8] $G\simeq L_3(q)$ has also only one conjugacy class of involutions, and the size
of that class is $\frac {|G|}{48}$. Hence in all these cases $i_2(G)<\frac {|G|}{20}$, which implies  by Lemma 3.2(2)
that $o(G)\geq 3.55$, as required.

It remains only to deal with $G\simeq L_3(q)$ for $q\geq 11$. We are going to use the information included in Lemma 2.13 in the paper [7]
and in the attached Table 2. If $G\simeq L_3(q)$ and $q\geq 11$, then $\frac {|G|}{20}>\frac 1{40}\frac {q^8}{q+1}$ and
$$i_2(G)\leq 2(q+1)q^4
\leq \frac {2(q+1)q^8}{(q^2-1)(q^2-1)}=\frac {2q^8}{(q-1)^2(q+1)}\leq \frac {2q^8}{100(q+1)}.$$
Hence $i_2(G)<\frac 1{40}\frac {q^8}{q+1}<\frac {|G|}{20}$
and by Lemma 3.2(2)
$o(G)\geq 3.55$, as required.
\qed
\enddemo

\proclaim{Lemma 4.4} Let $G\simeq U_3(q)$, with $q=3,4,5$ or $q\geq 8$ .
Then
$$o(G)\geq 3.55.$$
\endproclaim
\demo{Proof}
If $q\in \{3,4,8,9\}$, then by [8] $G\simeq U_3(q)$ has only one conjugacy class of involutions and the orders of their Sylow $2$-subgroups
are $\{2^5,2^6,2^9,2^5\}$, respectively.  If $q=5$, then by [8] $G\simeq U_3(5)$ has also only one conjugacy class of involutions, and the size
of that class is $\frac {|G|}{240}$. Hence in all these cases $i_2(G)<\frac {|G|}{20}$ and by Lemma 3.2(2)
$o(G)\geq 3.55$, as required.

It remains only to deal with $G\simeq U_3(q)$ for $q\geq 11$. By Lemma 2.13 and Table 2 in [7], $\frac {|G|}{20}>\frac 1{40}\frac {q^8}{q+1}$ and
$i_2(G)\leq 2(q+1)q^4$. Hence, as shown in the proof of Lemma 4.3, we have $i_2(G)<\frac {|G|}{20}$
and by Lemma 3.2(2)
$o(G)\geq 3.55$, as required.
\qed
\enddemo

We can now prove the following proposition.

\proclaim{Proposition  4.5} Let $G$ be a finite simple group with $|\Pi(G)|=3$ and suppose that $\Pi(G)\subseteq\{2,3,5,7, 11,13\}$. Then
$$ \text{either}\quad G\simeq A_5 \quad\text{or}\quad o(G)\geq 3.55.$$
\endproclaim
\demo{Proof} By Proposition 2.5(1), $G$ is isomorphic to one of the following groups:
$$A_5,A_6,S_4(3)\simeq U_4(2),L_2(7),L_2(8),L_3(3), U_3(3).$$
By Lemmas 4.1-4.4, if $G$ is isomorphic to one of the groups:
$$A_6,L_2(7),L_2(8),L_3(3), U_3(3)$$
then $o(G)\geq 3.55$, as required.

Finally, if $G\simeq U_4(2)$, then by [8]
$i_2(G)<2\frac {|G|}{96}<\frac{|G|}{20}$,
and by Lemma 3.2(2)
$o(G)\geq 3.55$, as required.
\qed
\enddemo

\proclaim{Lemma 4.6} Let either $G\simeq L_4(q)$ or $G\simeq U_4(q)$, with $q\geq 3$.
Then $o(G)\geq 3.55.$
\endproclaim
\demo{Proof} By  Lemma 2.13 and Table 2 in [7], $\frac {|G|}{20}>\frac 1{40}\frac {q^{15}}{q+1}$
and $i_2(G)\leq 2(q+1)q^8$. Since $q\geq 3$, we have
$$i_2(G)\leq 2(q+1)q^8 \leq \frac {2(q+1)q^{15}}{(q^2-1)(q^2-1)q^3}=\frac {2q^{15}}{(q-1)^2(q+1)q^3}\leq \frac {2q^{15}}{100(q+1)}.$$
Hence $i_2(G)<\frac 1{40}\frac {q^{15}}{q+1}<\frac {|G|}{20}$
and by Lemma 3.2(2)
$o(G)\geq 3.55$, as required.
\qed
\enddemo

\proclaim{Lemma 4.7} Let either $G\simeq L_5(q)$ or $G\simeq U_5(q)$, with $q\geq 2$.
Then $o(G)\geq 3.55.$
\endproclaim
\demo{Proof}  By  Lemma 2.13 and Table 2 in [7], $\frac {|G|}{20}>\frac 1{40}\frac {q^{24}}{q+1}$
and $i_2(G)\leq 2(q+1)q^{13}$. Since $q\geq 2$, we have
$$i_2(G)\leq 2(q+1)q^{13} \leq \frac {2(q+1)q^{24}}{(q^2-1)(q^2-1)q^7}=\frac {2q^{24}}{(q-1)^2(q+1)q^7}\leq \frac {2q^{24}}{100(q+1)}.$$
Hence $i_2(G)<\frac 1{40}\frac {q^{24}}{q+1}<\frac {|G|}{20}$
and by Lemma 3.2(2)
$o(G)\geq 3.55$, as required.
\qed
\enddemo

\proclaim{Lemma 4.8} Let  $G\simeq R(q)\simeq ^2G_2(q)$, with $q=3^{2n+1}$ and $n\geq 1$.
Then $o(G)\geq 3.55.$
\endproclaim
\demo{Proof} By  Lemma 2.13 and Table 2 in [7], $\frac {|G|}{20}>\frac 1{40}q^{7}$
and $i_2(G)\leq 2(q+1)q^3$. Since $q\geq 3^3=27$, we have
$$i_2(G)\leq 2(q+1)q^3 \leq \frac {2(q+1)q^7}{(q^2-1)q^2}=\frac {2q^7}{(q-1)q^2}\leq \frac {2q^7}{100}.$$
Hence $i_2(G)<\frac 1{40}q^7<\frac {|G|}{20}$
and by Lemma 3.2(2)
$o(G)\geq 3.55$, as required.
\qed
\enddemo

\proclaim{Lemma 4.9} Let  $G\simeq G_2(q)$, with  $q\geq 3$.
Then $o(G)\geq 3.55.$
\endproclaim
\demo{Proof} By  Lemma 2.13 and Table 2 in [7], $\frac {|G|}{20}>\frac 1{40}q^{14}$
and $i_2(G)\leq 2(q+1)q^7$. Since $q\geq 3$, we have
$$i_2(G)\leq 2(q+1)q^7 \leq \frac {2(q+1)q^{14}}{(q^2-1)q^5}=\frac {2q^{14}}{(q-1)q^5}\leq \frac {2q^{14}}{100}.$$
Hence $i_2(G)<\frac 1{40}q^{14}<\frac {|G|}{20}$
and by Lemma 3.2(2)
$o(G)\geq 3.55$, as required.
\qed
\enddemo

\proclaim{Lemma 4.10} Let  $G\simeq S_4(q)$, with $q\geq 5$.
Then $o(G)\geq 3.55.$
\endproclaim
\demo{Proof} By  Lemma 2.13 and Table 2 in [7], $\frac {|G|}{20}>\frac 1{80}q^{10}$
and $i_2(G)\leq 2(q+1)q^5$. Since $q\geq 5$, we have
$$i_2(G)\leq 2(q+1)q^5 \leq \frac {2(q+1)q^{10}}{(q^2-1)q^3}=\frac {2q^{10}}{(q-1)q^3}\leq \frac {2q^{10}}{200}.$$
Hence $i_2(G)<\frac 1{80}q^{10} <\frac {|G|}{20}$
and by Lemma 3.2(2)
$o(G)\geq 3.55$, as required.
\qed
\enddemo

\proclaim{Lemma 4.11} Let  $G\simeq Sz(2^{2m+1})\simeq \ ^2B_2(q)$,  $q = 2^{2m+1}$, with $m\geq 1$.
Then $o(G)\geq 3.55.$
\endproclaim
\demo{Proof} If $m=1$, then $G\simeq Sz(8)$ has only one conjugacy class of involutions and the order of its Sylow $2$-subgroups
is $2^6$. Hence  $i_2(G)<\frac {|G|}{20}$ and by Lemma 3.2(2)
$o(G)\geq 3.55$, as required.

So suppose that $m>1$ and $q\geq 32$.
By  Lemma 2.13 and Table 2 in [7], $\frac {|G|}{20}>\frac 1{40}q^{5}$
and $i_2(G)\leq 2(q+1)q^2$. Since $q\geq 32$, we have
$$i_2(G)\leq 2(q+1)q^2 \leq \frac {2(q+1)q^{5}}{(q^2-1)q}=\frac {2q^{5}}{(q-1)q}\leq \frac {2q^{5}}{100}.$$
Hence $i_2(G)<\frac 1{40}q^{5} <\frac {|G|}{20}$
and by Lemma 3.2(2)
$o(G)\geq 3.55$, as required.
\qed
\enddemo

\proclaim{Lemma 4.12} Let  $G\simeq S_6(q)$, with  $q\geq 2$.
Then $o(G)\geq 3.55.$
\endproclaim
\demo{Proof} By  Lemma 2.13 and Table 2 in [7], $\frac {|G|}{20}>\frac 1{80}q^{21}$
and $i_2(G)\leq 2(q+1)q^{11}$. Since $q\geq 2$, we have
$$i_2(G)\leq 2(q+1)q^{11} \leq \frac {2(q+1)q^{21}}{(q^2-1)q^8}=\frac {2q^{21}}{(q-1)q^8}\leq \frac {2q^{21}}{200}.$$
Hence $i_2(G)<\frac 1{80}q^{21}<\frac {|G|}{20}$
and by Lemma 3.2(2)
$o(G)\geq 3.55$, as required.
\qed
\enddemo

\proclaim{Lemma 4.13} Let  $G\simeq O_8^+(q)$, with  $q\geq 2$.
Then $o(G)\geq 3.55.$
\endproclaim
\demo{Proof} By  Lemma 2.13 and Table 2 in [7], $\frac {|G|}{20}>\frac 1{160}q^{28}$
and $i_2(G)\leq 2(q+1)q^{15}$. Since $q\geq 2$, we have
$$i_2(G)\leq 2(q+1)q^{15} \leq \frac {2(q+1)q^{28}}{(q^2-1)q^{11}}=\frac {2q^{28}}{(q-1)q^{11}}\leq \frac {2q^{28}}{400}.$$
Hence $i_2(G)<\frac 1{160}q^{28}<\frac {|G|}{20}$
and by Lemma 3.2(2)
$o(G)\geq 3.55$, as required.
\qed
\enddemo

\proclaim{Lemma 4.14} Let  $G\simeq \ ^3D_4(2)$.
Then $o(G)\geq 3.55.$
\endproclaim
\demo{Proof} By [8], $i_2(G)<2\frac {|G|}{3072}<\frac {|G|}{20}$ and by Lemma 3.2(2)
$o(G)\geq 3.55$, as required.
\qed
\enddemo

\proclaim{Lemma 4.15} Let  $G\simeq \ ^2F_4(2)'$.
Then $o(G)\geq 3.55.$
\endproclaim
\demo{Proof} By [8], $i_2(G)<2\frac {|G|}{1536}<\frac {|G|}{20}$ and by Lemma 3.2(2)
$o(G)\geq 3.55$, as required.
\qed
\enddemo

We can now prove the following propositions.

\proclaim{Proposition  4.16} Let $G$ be a finite simple group with $|\Pi(G)|=4$ and suppose that $\Pi(G)\subseteq\{2,3,5,7,11,13\}$. Then
$ o(G)\geq 3.55$.
\endproclaim
\demo{Proof} By Proposition 2.5(2), $G$ is isomorphic to one of the following groups:
$$ \align
L_2(q),\ &\text{$q$ a prime},\ |\pi(q^2-1)|=3,\\
L_2(2^m),\ &2^m-1=u,\ 2^m+1=3t,\ \text{$u,t$ are primes},\ t>3,\\
L_2(3^m),\ &3^m-1=2u,\ 3^m+1=4t,\ \text{$u,t$ are primes},\\
L_2(25),\ &L_2(49),\ L_3(4),\ L_4(3),\ U_3(4),\ U_3(5),\ U_4(3),\ U_5(2),\\
S_4(5),\ &
S_4(7),\ S_6(2),\ O_8^+(2),\ G_2(3),\ ^3D_4(2),\ ^2F_4(2)', Sz(8),\\
A_7,\ &A_8,\ A_9,\ A_{10}, M_{11},\ M_{12},\ J_2.
\endalign
$$
By the previous lemmas, only the following groups need to be checked: $M_{11}$, $M_{12}$ and $J_2$.

By [8], if $G\simeq M_{11}$, then  $i_2(G)=\frac {|G|}{48}<\frac{|G|}{20}$;
if $G\simeq M_{12}$, then $i_2(G)<2\cdot\frac{|G|}{192}< \frac {|G|}{20}$;
and if $G\simeq J_2$, then  $i_2(G)<2\cdot\frac{|G|}{240}< \frac {|G|}{20}$.

Hence by Lemma 3.2(2), each of these groups satisfies
$o(G)\geq 3.55$, as required.
\qed
\enddemo

\proclaim{Proposition  4.17} Let $G$ be a finite simple group with $|\Pi(G)|=5$ and suppose that $\Pi(G)\subseteq\{2,3,5,7,11,13\}$. Then
$ o(G)\geq 3.55$.
\endproclaim
\demo{Proof} By Proposition 2.5(3), $G$ is isomorphic to one of the following groups:
$$ \align
L_2(q),\ &\text{$q$ a prime power},\ |\pi(q^2-1)|=4,\\
L_3(q),\ &\text{$q$ a prime power },\ |\pi((q^2-1)(q^3-1))|=4,\\
U_3(q),\ &\text{$q$ a prime power },\ |\pi((q^2-1)(q^3+1))|=4,\\
O_5(q)&\simeq S_4(q),\ \text{$q$ a prime power },\ |\pi(q^4-1)|=4,\\
Sz(2^{2m+1})&\simeq \ ^2B_2(2^{2m+1}),\ \text{with}\ |\pi((2^{2m+1}-1)(2^{4m+2}+1))|=4,\\
R(q),\ &\text{with}\ q=3^{2m+1},\ |\pi(q^2-1)|=3\ \text{and}\ |\pi(q^2-q+1)|=1,\\
L_5(3),\ &S_6(3),\ U_4(5),\ U_6(2),\ O_7(3),\ O_8^+(3),\ G_2(4),\\
A_{11},\ &A_{12},\ M_{22},\ HS,\ McL.
\endalign
$$
By the previous lemmas, only the following groups need to be checked: $U_6(2)$, $O_7(3)$, $M_{22}$,  $HS$ and $McL$.

By [8], if $G\simeq U_6(2)$, then  $i_2(G)<3\cdot\frac{|G|}{36864}<\frac{|G|}{20}$;
if $G\simeq O_7(3)$, then  $i_2(G)<3\cdot\frac{|G|}{13824}<\frac{|G|}{20}$;
if $G\simeq M_{22}$, then  $i_2(G)=\frac {|G|}{384}<\frac{|G|}{20}$;
if $G\simeq HS$, then  $i_2(G)<2\cdot\frac{|G|}{2880}<\frac{|G|}{20}$;
and if $G\simeq McL$, then  $i_2(G)=\frac{|G|}{40320}<\frac{|G|}{20}$.

Hence by Lemma 3.2(2), each of these groups satisfies
$o(G)\geq 3.55$, as required.
\qed
\enddemo

\proclaim{Proposition  4.18} Let $G$ be a finite simple group with $|\Pi(G)|=6$ and suppose that $\Pi(G)\subseteq\{2,3,5,7,11,13\}$. Then
$ o(G)\geq 3.55$.
\endproclaim
\demo{Proof} By Proposition 2.5(4), $G$ is isomorphic to one of the following groups:
$$ \align
L_2(q),\ &\text{$q$ a prime power},\ |\pi(q^2-1)|=5,\\
L_3(q),\ &\text{$q$ a prime power },\ |\pi((q^2-1)(q^3-1))|=5,\\
L_4(q),\ &\text{$q$ a prime power },\ |\pi((q^2-1)(q^3-1)(q^4-1))|=5,\\
U_3(q),\ &\text{$q$ a prime power },\ |\pi((q^2-1)(q^3+1))|=5,\\
U_4(q),\ &\text{$q$ a prime power },\ |\pi((q^2-1)(q^3+1)(q^4-1))|=5,\\
O_5(q)&\simeq S_4(q),\ \text{$q$ a prime power },\ |\pi(q^4-1)|=5,\\
G_2(q),\ &\text{$q$ a prime power },\ |\pi(q^6-1)|=5,\\
Sz(2^{2m+1})&\simeq \ ^2B_2(2^{2m+1}),\ \text{with}\ |\pi((2^{2m+1}-1)(2^{4m+2}+1))|=5,\\
R(3^{2m+1})&,\ \text{with}\ |\pi((3^{2m+1}-1)(3^{6m+3}+1))|=5,\\
L_6(3),\ &A_{13},\ A_{14},\ A_{15},\ A_{16,},\ Suz,\ Fi_{22}.
\endalign
$$
By the previous lemmas, only the following groups need to be checked: $L_6(3)$, $Suz$ and $Fi_{22}$.

If $G\simeq L_6(3)$,  then by  Lemma 2.13 and Table 2 in [7] $\frac {|G|}{20}>\frac 1{160}3^{35}>3^{21}$
and
$i_2(G)< 8\cdot3^{19}<3^{21}<\frac {|G|}{20}$.
Hence, by Lemma 3.2(2), $o(G)\geq 3.55$, as required.

If $G\simeq Suz$, then by [8] $i_2(G)<2\cdot \frac {|G|}{161280}<\frac {|G|}{20}$ and by
Lemma 3.2(2) $o(G)\geq 3.55$, as required.

Finally, if $G\simeq Fi_{22}$, then by [8] $i_2(G)<3\cdot \frac {|G|}{1769472}<\frac {|G|}{20}$ and by
 Lemma 3.2(2) $o(G)\geq 3.55$, as required.
\qed
\enddemo

We close this section with the following concluding theorem.
\proclaim {Theorem 4.19} Let $G$ be a finite simple group with $\Pi(G)\subseteq \{2,3,5,7,11,13\}$.
Then either $G\simeq A_5$ or $o(G)\geq 3.55> o(A_5)$.
\endproclaim
\demo{Proof} Since $G$ is a simple group, it follows that $|\Pi(G)|\in \{3,4,5,6\}$. Therefore by
Propositions 4.5, 4.16, 4.17 and 4.18, either $G\simeq A_5$ or $o(G)\geq 3.55> o(A_5)$, as required.
\qed
\enddemo

\heading 5. Non-solvable groups with an automorphism inverting many elements.\\
\endheading

In their paper [32], H. Liebeck and D. MacHale  proved that if a finite group $G$ admits an automorphism
which inverts more than one half of the elements of $G$ (i.e. sends more than one half of its elements to
their inverses), then $G$ is solvable. Later, W.M. Potter proved in [35] that  if a finite non-solvable group $G$
admits an automorphism
which inverts more than one quarter of the elements of $G$, then $G\simeq B\times A_5$, where $B$ denotes a finite
abelian subgroup of $G$. In this section, we generalize Potter's theorem and, using his methods and his ideas,
we prove the following theorem.
\proclaim{ Theorem 5.5} Let $G$ be a finite non-solvable group,  which contains no non-trivial solvable normal subgroups and let
$\theta \in Aut(G)$. If $\theta$ inverts more than $\frac 29|G|$ of the elements of $G$, then $G\simeq A_5$.
\endproclaim

We shall use the following notation from [35]. Let $G$ be a finite group and let $\theta \in Aut(G)$.
We shall write
$$S(\theta)=\{g\in G\mid g^{\theta}=g^{-1}\}\quad \text{and}\quad r(G,\theta)=\frac {|S(\theta)|}{|G|}.$$
When $x\in G$, we shall denote by $x\theta$ the automorphism of $G$ defined by
$$g^{x\theta}=(x^{-1}gx)^{\theta}\quad \text{for all}\quad g\in G.$$

We shall use also the following four results from [35].

\proclaim {Lemma 5.1} Let $G$ be a finite group, $\theta \in Aut(G)$ and let $H$ be a subgroup of $G$. Then
the following statements hold:

(i) if $x\in S(\theta)$, then $S(x\theta)=S(\theta)x^{-1}$,

(ii) there exists $x\in S(\theta)$ so that
$$|Hx\cap S(\theta)|=|H\cap S(x\theta)|\geq r(G,\theta)\cdot |H|.$$
\endproclaim
\demo{Proof} See Lemma 2.1 of [35].
\qed
\enddemo

If $\theta \in Aut(G)$ and $N$ is a normal $\theta$-invariant subgroup of $G$, then we may define an
automorphism $\bar\theta$ of $G/N$ by $(Ng)^{\bar\theta}=Ng^{\theta}$ for all $Ng\in G/N$.
The following lemma provides a useful relationship between $r(G,\theta)$ and $r(G/N,\bar\theta)$.

\proclaim {Lemma 5.2} Let $N$ be a normal subgroup of a finite group $G$ and let  $\theta \in Aut(G)$ satisfy
$N^{\theta}=N$. If for all $g\in S(\theta)$ we have $|Ng\cap S(\theta)|\leq t|N|$ for some $0<t\leq 1$, then
$$t^{-1}\cdot r(G,\theta)\leq r(G/N,\bar\theta).$$
\endproclaim
\demo{Proof} See Lemma 2.2 of [35].
\qed
\enddemo

\proclaim {Lemma 5.3} Let $H$ be a subgroup of a finite group $G$ and let $p$ be the smallest prime divisor of $|H|$.
If $\theta\in Aut(G)$ inverts $H$ pointwise and if $r(G,\theta)>\frac 1p$, then
$$|G:C_G(H)|\leq \frac {p-1}{p\cdot r(G,\theta)-1}.$$
\endproclaim
\demo{Proof} See Theorem 2.3 of [35].
\qed
\enddemo

\proclaim {Lemma 5.4} Let $G$ be a finite group, $\theta \in Aut(G)$ and let $P$ be a Sylow $p$-subgroup of $G$.
If $r(G,\theta)>\frac 2{p+1}$, then $P$ is normal in $G$.
\endproclaim
\demo{Proof} See Theorem 2.5 of [35].
\qed
\enddemo

We are now ready for the proof of Theorem 5.5.

\demo{Proof of Theorem 5.5} Suppose that $\theta \in Aut(G)$ and $r(G,\theta)>\frac 29$.

Suppose first that $G$ is a  simple group. If $p$ is a prime dividing $|G|$  and $P$ is a Sylow $p$-subgroup of $G$,
then $P$ is not normal in $G$ and by Lemma 5.4
$$\frac 2{p+1}\geq r(G,\theta)>\frac 29.$$
Hence $p+1<9$ and $p\in \{2,3,5,7\}$. Since $G$ is non-solvable, we have $|\Pi(G)|\geq 3$.
Thus $G$ has a non-trivial Sylow $p$-subgroup for either $p=5$ or $p=7$. Since  $p\geq 5$,
it follows by Lemma 5.1(2) that we can choose $\theta$ such that
$$|P\cap S(\theta)|> \frac 29|P|>\frac{|P|}p.$$
Thus $\langle P\cap S(\theta) \rangle =P$,
which implies that $P^{\theta}=P$ and the restriction of $\theta$ to $P$ has order $2$. As $p$ is odd, it follows by Lemma 4.1 in Chapter 10 of [10]
that $P=C_P(\theta)(S(\theta)\cap P)$. Hence $|C_P(\theta)|=\frac {|P|}{|S(\theta)\cap P|}<p$, so $S(\theta)\cap P=P.$
Therefore $P$  is inverted pointwise by $\theta$ and hence it is abelian. Since $r(G,\theta)>\frac 29>\frac 1p$, it follows by Lemma 5.3 that
$$|G:C_G(P)|\leq \frac {p-1}{p\cdot r(G,\theta)-1}<\frac {p-1}{p\cdot \frac 29-1}=\frac {9(p-1)}{2p-9}.$$
If $p=5$, then $|G:C_G(P)|< 36$ and if $p=7$, then $|G:C_G(P)|< \frac {54}5<11$. Since $P$ is abelian and $G$ is not
$p$-nilpotent, it follows by Burnside's Lemma that
$P\leq C_G(P)<N_G(P)$. If $p=7$, then $|G:N_G(P)|\geq 8$ and $|G:C_G(P)|\geq 16$, a contradiction.

Therefore $p=5$,  $\Pi(G)=\{2,3,5\}$ and by Proposition 2.5(1) $G$ is isomorphic to one of the groups $ \{A_5,A_6, S_4(3)\}$. By [8], in
each case $|P|=|C_G(P)|=5$  and if
$G$ is isomorphic to either $A_6$ or $S_4(3)$, then $|G:C_G(P)| > 36 $, a contradiction. Hence $G\simeq A_5$, as required.

Now suppose that $G$ is  a non-simple group and let $N$ be a minimal normal subgroup of $G$.
Since $G$ contains no non-trivial solvable normal subgroups, $N$ is a direct product of isomorphic  simple groups. It follows from Lemma 5.1(2) that we
can choose an automorphism $\theta$  of $G$ such that $|N\cap S(\theta)|> \frac 29|N|>\frac 15|N|$.
Since a perfect group cannot have a proper subgroup of index $\leq 4$, it follows that $N=\langle N\cap S(\theta)\rangle$ and $N^{\theta}=N$.
Let $X$ be one of the simple factors of $N$. Arguing similarly on $X$, we may conclude that  $X^{\theta}=X$ and  $|X\cap S(\theta)|> \frac 29|X|.$
It follows then, by the first part of the proof, that $X\simeq A_5$.

Since $X$ is non-solvable, it follows by Lemma 5.1(1) and Proposition 2.3(1) that
$$|X\cap S(y\theta)|=|Xy\cap S(\theta)|\leq \frac 4{15}|X|\quad \text{for every} \quad  y\in N\cap S(\theta).$$
Thus, by Lemma 5.2,
$$r(N/X,\bar\theta)> \frac {15}4\cdot\frac 29=\frac 56>\frac 4{15},$$
which implies by Proposition 2.3(1) that $N/X$ is solvable. Hence $N=X\simeq A_5$.
Since $X$ is normal in $G$, it follows as before that $G/X$ is solvable. Since $X$ is simple, we have $X\cap C_G(X)=1$, which implies
that $C_G(X)$ a solvable normal subgroup of $G$. But
$G$ contains no non-trivial solvable normal subgroups, so
$C_G(X)=\{1\}$. Thus $G$ is isomorphic to a subgroup of $Aut(A_5)\simeq S_5$.
But $S_5$ is a complete group and it has no inner
automorphisms inverting more than $\frac 16$ of its elements. Since $\frac 16<\frac 29$, it follows that $X\leq G$ and $|G|<|S_5|$.  Hence $G=X\simeq A_5$,
in  contradiction to our assumption that $G$ is not a simple group.

The proof of Theorem 5.5 is now complete.
\qed
\enddemo

\heading 6. Proofs of Theorem A and Theorem B.\\
\endheading

We start this section with the following lemma.
\proclaim {Lemma 6.1} Let $G$ be a finite non-solvable group and suppose that there exists a
normal subgroup $N$ of $G$ with $|G/N|=3$. Then
$$o(N)<3(o(G)-2.66).$$
In particular
$$o(G)> 3.66.$$
\endproclaim
\demo{Proof} Write $Y=G\setminus N$. Then $|Y|=\frac 23|G|$ and for every $y\in Y$ we have $y^3\in N$. Since $y\notin N$,
$|yN|=3$ in $G/N$ and by Lemma 3.1(1)  $|y|=3h$
for some positive integer $h$. Write $T=\{y\in Y\mid |y|=3\}$ and $V=Y\setminus T$. Then $\psi(Y)\geq 3|T|+6|V|$. Moreover, $i_3(G)\geq |T|$ and
$i_3(G)<\frac 7{20}|G|$ by Proposition 2.4. Suppose that $|T|>\frac 23|Y|$. Then
$$\frac 7{20}|G|>i_3(G)\geq |T|>\frac 23|Y|=\frac 23\cdot \frac 23|G|=\frac 49|G|,$$
which is a contradiction, since  $\frac 7{20}<\frac 49$. Therefore $|T|\leq \frac 23|Y|$ and hence $|V|=|Y|-|T|\geq |Y|-\frac 23|Y|=\frac 13|Y|$.
It follows that
$$\psi(G)-\psi(N)=\psi(Y)\geq 3|T|+6|V|=3(|T|+|V|)+3|V|\geq 3|Y|+|Y|=\frac 83|G|,$$
which implies that
$$o(G)=\frac {\psi(G)}{|G|}\geq \frac {\psi(N)}{|G|}+\frac 83> \frac 13o(N)+2.66.$$
Thus
$$o(N)< 3(o(G)-2.66)$$
as required.
Since $N$ is non-solvable, it follows by Lemma 3.2(1) that $o(N)\geq 3.11$. Hence $3<o(N)<3(o(G)-2.66)$, which
implies that $1<o(G)-2.66$. So
$o(G)>3.66$,
as required.
\qed
\enddemo

We proceed now with the proofs of Theorem A and Theorem B.
\demo{Proof of Theorem A} Suppose that $o(G)\leq \frac {13}6=o(S_3)$. Then by Lemma 3.2(1) $G$ is a solvable group, and hence $G$
contains a normal subgroup $N$ with $|G/N|=p$, where $p$ is a prime. If $|N|=1$, then $|G|=p$ and by our assumptions
$o(G)=\psi(G)/|G|=\frac {p(p-1)+1}p=p-1+\frac 1p\leq \frac {13}6$. Thus $p=2$,  $|G|=2$ and $o(G)=\frac 32<\frac {13}6$, as required.
So assume that $|N|>1$.

We claim that $p=2$ and $G=N\dot\cup xN$, where $x\in G\setminus N$ and $|x|=2$.

Indeed, since $\psi(G/N)=p(p-1)+1$, we have $o(G/N)=\frac {p(p-1)+1}p=p-1 +\frac 1p$ and by Lemma 3.1(2) $p-1+\frac 1p<o(G)\leq \frac {13}6$.
Hence $p=2$ and $G=N\dot\cup xN$, with $x\in G\setminus N$. By Lemma 3.1(1), $2$ divides $|xn|$ for each $n\in N$.
If $|xn|\geq 4$ for each $n\in N$, then $\psi(G)\geq \psi(N) +4|N|$ and $o(G)\geq \frac 12o(N)+2$.
Hence $o(N)\leq 2(o(G)-2)\leq 2(\frac  {13}6-2)=\frac 13<1$, a contradiction. It follows that there exists $n\in N$  with $|xn|=2$,
and hence we may assume that $|x|=2$. The proof of the claim is now complete.

Next we claim that $N$ is abelian. Recall that $2$ divides $|xn|$ for each $n\in N$. This implies that $\psi(G)\geq \psi(N)+2|N|$.
If $N$ is non-abelian, then by Proposition 2.3(2) the number of elements $n$ of $N$ which are inverted by $x$, i.e. such that
$|xn|=2$, is at most $\frac 34|N|$. Write $V=\{n\in N\mid |xn|\geq 4\}$. Then $|V|\geq (1-\frac 34)|N|=\frac 14|N|$ and
$$\psi(G)\geq \psi(N)+2(|N|-|V|)+4|V|=\psi(N)+2|N|+2|V|\geq \psi(N)+2|N|+2(\frac 14|N|).$$
Hence $o(G)\geq \frac 12 o(N)+1+\frac 14=\frac 12 o(N)+\frac 54$ and
$$o(N)\leq 2(o(G)-\frac 54)\leq 2(\frac {13}6-\frac 54)=2(\frac {11}{12})<2.$$
But then $N$ is an elementary abelian $2$-group, in contradiction to our assumption. Hence $N$ is an abelian group, as claimed.

Next we claim that $o(N)\leq \frac 73$ and  $N$ is either an abelian $2$-group or $|N|=3$.

Indeed, since $\psi(G)\geq \psi(N)+2|N|$, we have $o(G)\geq \frac 12 o(N)+1$ and hence
$$o(N)\leq 2(o(G)-1)\leq 2(\frac {13}6-1)=\frac 73,$$
as claimed.

Suppose that $N=D\times O$, where $D\neq 1$ is a $2$-group and $O\neq 1$ is a group of odd order.
Then
$o(N)\geq \frac 72$ by Lemma 1.1(3),  in contradiction to $o(N)\leq\frac 73$. Hence $N$ is either an abelian $2$-group
or a group of odd order. In the second case $\frac 73\geq o(N)\geq 3-\frac 2{|N|}$ by Lemma 1.1(2), which implies that $\frac 2{|N|}\geq \frac 23$
and $|N|\leq 3$. Hence  $N$ is either an abelian  $2$-group or $|N|=3$, as claimed.

First suppose that $N$ is an abelian $2$-group. We claim that $N$ is an elementary abelian $2$-group. Indeed, if this is not the case, then
$N$ contains a normal subgroup $K$ such that $N/K\simeq C_4$ and $o(N/K)=\frac {11}4$. If $|K|=1$, then $o(N)= \frac {11}4>\frac 73$,
in contradiction to $o(N) \leq \frac 73$. If $|K|>1$, then by Lemma 3.1(2) $\frac {11}4=o(N/K)<o(N)\leq \frac 73$, a contradiction. Hence $N$ is an
elementary abelian $2$-group, as claimed.

We now claim that also $G$ is an elementary abelian $2$-group. Let
$W=\{n\in N\mid |nx|=2\}=C_N(x)$. If $W<N$, then $|W|\leq \frac 12|N|$ and hence $|N\setminus W|\geq \frac 12|N|$.
It follows that
$$\psi(G)\geq \psi(N)+2|W|+4|N\setminus W|= \psi(N)+2|N|+2|N\setminus W|
\geq \psi(N)+3|N|,$$
which implies that
$o(N)\leq 2o(G)-3\leq 2\cdot\frac {13}6-3=\frac 43$, in contradiction to Lemma 1.1(1).
Therefore $C_N(x)=W=N$ and hence $G$ is an elementary abelian
$2$-group, satisfying $o(G)<2<\frac {13}6$, as required.

Finally suppose that $|N|=3$. Then $|G|=6$ and since $o(C_6)=\frac 72>\frac {13}6$, it follows that $G\simeq S_3$ and $o(G)=\frac {13}6$.

The proof of Theorem A is now complete.
\qed
\enddemo

We conclude this paper with the proof of Theorem B.

\demo{Proof of Theorem B} Suppose that $G$ is a non-solvable finite group,
which satisfies $o(G)\leq o(A_5)$ and is not isomorphic to $A_5$, and
suppose that $G$ is of minimal order satisfying these conditions.
We shall reach a contradiction, which will indicate that if a finite group $G$
satisfies $o(G)\leq o(A_5)$, then either $G\simeq A_5$ or $G$ is solvable, as required.

First we show that $G$ is a non-simple group.

Suppose that $G$ is a simple group. First we restrict the number of primes dividing $|G|$. Let $p$ be a prime dividing $|G|$.  By Lemma 3.3(1),
if $p\geq 17$ and $G$ is a non-solvable group satisfying $o(G)\leq o(A_5)$, then $G$ is $p$-solvable. Since $G$ is a simple group
satisfying $o(G)\leq o(A_5)$, it follows
that $\Pi(G)\subseteq \{2,3,5,7,11,13\}$  and by Theorem 4.19 either $G\simeq A_5$ or $o(G)\geq 3.55>o(A_5)$, in contradiction to our assumptions.

Hence $G$ is a non-simple non-solvable group and there exists a normal subgroup $S$ of $G$, satisfying $1<|S|<|G|$. By Lemma 3.1(2) we have $o(G/S)<o(G)\leq o(A_5)$,
so $G/S$ is not isomorphic to $A_5$ and since $o(G/S)<o(A_5)$, it follows from the minimality of $G$ that $G/S$ is solvable. Therefore there exist
a normal subgroup $N$ of $G$ and a prime $p$, such that   $|G/N|=p$ and $o(G/N)<o(G)\leq o(A_5)=3.51666\dots<3.52$. Since $G$ is non-solvable,
also $N$ is non-solvable.

If $p>3$, then
$$o(G/N)=o(C_p)=\frac {p(p-1)+1}p=(p-1)+\frac 1p>4>o(A_5),$$
a contradiction. Hence either $p=2$ or $p=3$. If $p=3$, then by Lemma 6.1 we have $o(G)\geq 3.66$, in contradiction to our assumption that
$o(G)<3.52$.

Hence $p=2$ and $|G/N|=2$. Thus $G=N\dot\cup xN$ for some $x\in G\setminus N$, with $|xN|=2$ in $G/N$. Therefore
$$\psi(G)=\psi(N)+\psi(xN)$$
and by Lemma 3.1(1), $2$ divides $|xn|$ for each $n\in N$. If $|xn|\geq 4$ for each $n\in N$, then $\psi(xN)\geq 4|N|$
and $\psi(G)\geq \psi(N)+4|N|$. Therefore  $o(G)=\frac {\psi(G)}{|G|}\geq \frac {\psi(N)}{|G|}+2=\frac 12o(N)+2$ and
$o(N)\leq 2(o(G)-2)<2\cdot (3.52-2)=3.04$, in contradiction to Lemma 3.2(1). Hence there exists $xn\in xN$ of order $2$ and we may assume that $|x|=2$.

Write $X=\{xn\mid n\in N, |xn|=2\}$.  Obviously $|xn|=2$ if and only if $n^x=n^{-1}$ and $|X|=|\{n\in N\mid n^x=n^{-1}\}|$. Moreover,
$$\psi(xN)\geq 2|X|+4(|xN\setminus X|).\tag{$*$}$$
If $N$ contains a non-trivial normal solvable subgroup $K$,
then $M=KK^x$ is a non-trivial normal solvable subgroup of $G$ and by Lemma 3.1(2) we have $o(G/M)<o(G)\leq o(A_5)$.
Since $G/M$ is a non-solvable, it follows by the minimality of $G$ that $G/M\simeq A_5$, a contradiction. Hence $N$ contains no non-trivial
solvable normal subgroups and $N\cap C_G(N)=1$. Moreover, if $\theta$ denotes the conjugation of $N$ by $x$, then Theorem 5.5 implies that
either $N\simeq A_5$ or $|X|\leq \frac 29|N|$.

First assume that $N\simeq A_5$. Since $N\cap C_G(N)=1$, if $C_G(N)\neq 1$, then $|C_G(N)|=2$ and $G=N\times C_G(N).$  Hence, by Lemma 3.1(2),
we get $o(N)=o(G/C_G(N))<o(G)\leq o(A_5)$, a contradiction, since $N\simeq A_5$. Therefore $C_G(N)=1$ and $G=G/C_G(N)\leq Aut(A_5)=S_5$. Since $|G|=2|N|=|S_5|$,
it follows that
$G=S_5$. But $\psi(S_5)=501$ and $|S_5|=120$, so $o(S_5)=\frac {501}{120}=4.175>o(A_5)$,  in contradiction to our assumption that $o(G)\leq o(A_5)$.

Finally, assume that $|X|\leq \frac 29|N|$. Then $|xN\setminus X|\geq |xN|-\frac 29|N|=\frac 79|N|$ and it follows by ($*$) that
$$\psi(G)=\psi(N)+\psi(xN)\geq \psi(N)+2|X|+4(|xN\setminus X|)=\psi(N)+2(|xN\setminus X|) +2|N|.$$
Hence
$\psi(G)\geq \psi(N)+2|N|+2\cdot \frac 79|N|$ and
$$o(G)\geq\frac 12 o(N)+1+\frac 79> \frac 12o(N)+1.777.$$
It follows that $N$ is a non-solvable group satisfying
$$o(N)< 2(o(G)-1.777)<2(3.52-1.777)=2\cdot 1.743=3.486<o(A_5).$$
Since $N<G$ and it is not isomorphic to $A_5$, we obtained
a contradiction to the minimality of $G$.

The proof of Theorem B is now complete.
\qed
\enddemo

\enddocument

\enddocument